\documentclass[a4paper,12pt]{amsart}
  \usepackage{macros_anglais}

\title[The Schoen--Webster Theorem]{Global rigidity in CR geometry: 
the Schoen--Webster Theorem}

\newcommand{\ssp}{{\mathcal S}}

\newcommand{\Heis}{\mathcal H}
\DeclareMathOperator{\T}{Tor}
\DeclareMathOperator{\Rm}{Rm}
\DeclareMathOperator{\R}{R}
\newcommand{\restric}[1]{\vert_{#1}}
\DeclareMathOperator{\LL}{\mathcal{L}}
\DeclareMathOperator{\vol}{vol}

\begin{document}

\maketitle

\begin{center}
Beno\^{\i}t KLOECKNER  \footnote{bkloeckn@fourier.ujf-grenoble.fr,
Institut Fourier, 100 rue des Maths, BP 74, 38402 St Martin d'Hy\`eres, France.} 
and Vincent MINERBE \footnote{minerbe-v@univ-nantes.fr, Laboratoire
 de math\'ematiques Jean Leray,
Universit\'e de Nantes, 2 rue de la Houssini\`ere, BP 92208,
44322 Nantes, cedex 3, France.}
\end{center}

\begin{abstract}
Schoen-Webster theorem asserts a pseudoconvex CR manifold whose automorphism
group acts non properly is either the standard sphere or the Heisenberg space.
The purpose of this paper is to survey successive works around this result and 
then provide a short geometric proof in the compact case.
\end{abstract}

\vskip 0.5cm


\textit{MS classification numbers}: 32V05, 32V20, 53C24.

\vskip 0.5cm

Among the many aspects of geometric rigidity, the vague principle
according to which a given geometry is \defini{rigid} when
``few manifolds admit a large automorphism group'' has a fairly rich history. 
In this survey paper, we try to show how 
strictly pseudoconvex CR geometry fits into this concept of
rigidity.

Andr\'{e} Lichnerowicz first raised the question in the conformal case. 
It is well known that the isometry group of a compact Riemannian manifold 
is compact, due to the compactness of the
group $\On{n}$ (see Section \ref{forigidity}). Since the corresponding group 
$\CO{n}$ of conformal geometry is not compact, one might expect some compact 
manifolds to have noncompact conformal groups. There is a simple example:
the Euclidean sphere has conformal group $\SO{1,n+1}$. The Lichnerowicz
conjecture stating that there are no other examples was settled in the early
seventies by Jacqueline Ferrand \cite{Ferrand1} and in a weak form
by Morio Obata\footnote{It appeared later that the proof was flawed
at some point, but Jacques Lafontaine gave a corrected proof
in  \cite{Lafontaine}} \cite{Obata}; it was extended by Ferrand
a while later \cite{Ferrand2}.

A few years after Obata and Ferrand's works, it appeared
that Lichnerowicz conjecture was not specific to conformal geometry:
Sidney Webster extended parts of the proof of Obata
in the setting of (strictly pseudoconvex) CR geometry \cite{Webster1}. The question
raised a lot of interest again in the nineties, several mathematicians
trying to work their way out from Webster's result to the full statement.
Richard Schoen gave the first complete proof, using original 
analytic methods related to the Yamabe problem \cite{Schoen}. In fact, he gave a proof in the 
conformal case that adapts to CR geometry and obtained the following result.\footnote{We
chose to name it after both Webster, who initiated the topic, and Schoen,
who gave the first complete proof.}

\begin{theo*}[Schoen--Webster]
Let $M$ be a strictly pseudoconvex CR manifold, not necessarily compact. 
If its automorphism group $\Aut(M)$ acts non-properly, then $M$ is either 
the standard CR sphere $\ssp$ or $\ssp$ with one point deleted.
\end{theo*}

Let us recall that an action  of a topological group $G$ is \defini{proper}
if for any compact subset $K$ of $M$, the subset
\[G_K=\ensemble{g\in G}{g(K)\cap K\neq\varnothing}\]
of $G$ is compact. In particular if $M$ is compact, $\Aut(M)$ acts properly if and
only if it is compact.

The paper is organised as follows. The first section is devoted
to preliminaries, including CR geometry, two properties that are important in the
sequel and $(G,X)$-structures. We then survey the successive works on the Schoen--Webster 
Theorem, 
trying to give for (almost) each result the flavor of the proof without getting
into too much detail. The word ``proof'' will therefore often be followed
by quite imprecise arguments.
The last section is devoted to a cleaned geometric proof of the theorem 
when $M$ is compact, based on some of the ideas exposed.

Before getting started, let us point out that Bun Wong proved  
a very close theorem for domains of $\mC^{n+1}$ \cite{Wong}.
Many developments arose from his result and parts of the Schoen--Webster Theorem can
be deduced from this work. Indeed, unless $n=1$, a compact strictly pseudoconvex
CR manifold $M^{2n+1}$ can always be embedded as the boundary of a domain of 
$\mC^{n+1}$
and its automorphisms can be extended to automorphisms of the domain. See
\cite{Lee} for details due to Daniel Burns.

However, we will not
discuss Wong's theorem and its improvements. First, it cannot be of any help for
the least dimensional case. Second, we are interested
in more intrinsic methods of proof, independant of any embedding. For further 
informations on this topic, the reader should refer to \cite{WongSurvey}.

For the sake of completeness, note that in \cite{Pansu} Pierre Pansu 
gave a hint of how one could try to adapt Ferrand's proof to the CR case.

\section{Preliminaries}

\subsection{Basics of CR geometry}

We only give a glimpse on CR geometry. The interested reader can refer to
\cite{Dangelo} or \cite{Jacobowitz}.

Given a $2n+1$-dimensional manifold $M$, a \defini{CR structure} on $M$ is a couple 
$(\xi,J)$ where:
\begin{enumerate}
\item $\xi$ is a $2n$-dimensional subbundle of $TM$,
\item $J$ is a pseudocomplex operator on $\xi$: 
      \[J_x : \xi_x \to \xi_x, \qquad J_x^2 = -\id \quad\forall x\in M,\]
\item\label{integrability} for all vector fields $X$, $Y$ tangent to $\xi$, the vector
      field $[JX,Y]+[X,JY]$ is tangent to $\xi$ and the following \defini{integrability
      condition} holds:
      \begin{equation}
        J([JX,Y]+[X,JY])=[JX,JY]-[X,Y].\nonumber
      \end{equation}
\end{enumerate}

Any smooth hypersurface $H$ in a complex manifold $X$ admits a natural CR-structure:
denoting by $J$ the complex structure of $X$, one can define $\xi=TH\cap J(TH)$
so that $J$ acts on $\xi$; note that the vanishing of the
Nijenhuis tensor implies the integrability condition.

A differentiable map between two CR manifolds is a \defini{CR map} if it conjugates
the hyperplanes distributions and the pseudocomplex operators.
An \defini{automorphism} of a CR manifold $M$ is a diffeomorphism of $M$
that is a CR map. The group of those is denoted by $\Aut(M)$
and its identity component by $\Aut_0(M)$.

\subsubsection{Calibrations, the Levi form and the Webster metric}

Given a CR structure $(\xi,J)$ on $M$, a (possibly local) $1$-form $\theta$ 
such that $\xi=\ker \theta$ is called a (local) \defini{calibration}.
One can always find local calibrations. If
$M$ is orientable, one can
always find a global calibration. However, a calibration need not be 
preserved by automorphisms.

\emph{From now on, all the manifolds under consideration are assumed to
be connected and orientable; all the calibrations are assumed to be global.}

Given a calibration $\theta$, one defines on $\xi$ the \defini{Levi form}:
\[L_\theta(\cdot)=d\theta(\cdot,J\cdot).\]
As a consequence of the integrability condition, the Levi form is a
quadratic form.

A change of calibration induces a linear change in the Levi form:
\begin{equation}
  L_{\lambda\theta}=\lambda L_\theta,
\end{equation}
thus its signature is, up to a change of sign, a CR invariant.

A CR structure is said to be \defini{strictly pseudoconvex} if
its Levi form is definite (and then we choose our calibrations so that
it is positive definite). It implies that $d\theta$ is nondegenerate
on $\ker\theta$, that is $\theta$ is a contact form.
If the Levi form vanishes at each point, the CR structure is said
to be \defini{Levi-flat}. Then $d\theta$ is zero on $\xi$ and the
Frobenius Theorem shows that $\xi$ defines a foliation.

One should therefore not think of CR geometry as one geometry : each signature
of the Levi form corresponds to a geometry of its own, just like Lorentzian
and Riemannian geometry (or foliations and contact structures)
are related, but different kind of geometries.

\begin{figure}\begin{center}
\input{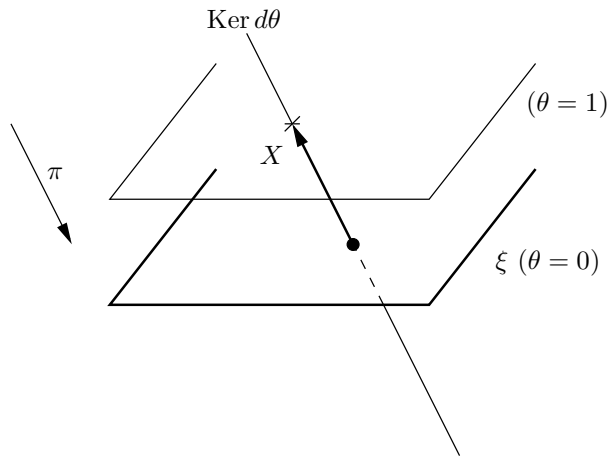}
\caption{The Reeb vector field of a calibration}\label{champ_reeb}
\end{center}\end{figure}

Given a calibration $\theta$ on a strictly pseudoconvex CR manifold, there is a 
single vector field $X$, called the \defini{Reeb vector field} of $\theta$, 
that satisfies:
\begin{equation}
\theta(X) = 1\quad\text{and}\quad X\contraction d\theta=0.
\end{equation}
See Figure \ref{champ_reeb}.

Denote by $\pi:TM\rightarrow\xi$ the linear projection on
$\xi$ along the direction of $X$. If the Levi form is positive definite,
one gets a Riemannian metric on $M$ called the \defini{Webster metric}:
\begin{equation}
W_\theta=L_\theta\rond\pi+\theta^2.
\end{equation}
A change of calibration $\theta'=\lambda\theta$ changes the metric
by a factor $\lambda$ along $\xi$ and by a factor $\lambda^2$ 
``transversally'' that is, on the quotient $TM/\xi$. Therefore,
the Webster metric does not define a canonical conformal structure
on a CR manifold.

Note that if $M$ has dimension $2n+1$, the calibration $\theta$ defines a volume form $\theta \wedge d\theta^n$ 
which is compatible with the Webster metric. 

\subsubsection{The Webster scalar curvature and the pseudoconformal Laplacian}

There is also a natural metric connection $\nabla_\theta$ on $TM$, the so called Tanaka-Webster connection; beware its torsion 
$\T_\theta$ does not vanish in general. Contracting the curvature $\Rm_\theta$ of this connection along $\xi$, we obtain a scalar 
curvature $\R_\theta$. A subelliptic Laplacian $\Delta_\theta$ arises by taking (minus) the trace over $\xi$ of the Hessian 
corresponding to $\nabla_\theta$; the following integration by parts formula holds:
$$
\forall u,v\in C^\infty_c(M),\quad \int_M (\Delta_\theta u) v \theta \wedge d\theta^n
= \int_M L_\theta \left(du\restric{\xi},dv\restric{\xi} \right) \theta \wedge d\theta^n,
$$
where $\xi$ and $\xi^*$ are identified thanks to $L_\theta$. To understand the relevance of this operator, consider another calibration 
$\theta'$, which we write $\theta' = u^{\frac{2}{n}} \theta$ for some smooth positive function $u$. The scalar curvature then transforms 
according to the following law: 
$$
R_{\theta'} = b(n)^{-1} u^{-\frac{n+2}{n}} \LL_{\theta} u
$$
where $b(n)=\frac{n+1}{4n+2}$ and $\LL_{\theta} = \Delta_{\theta} + b(n) \R_{\theta}$. 

This formula is pretty similar to a conformal one. Indeed, given conformally equivalent metrics $g$ and $h= u^{\frac{4}{n-2}} g$, 
for some positive function $u$, the Riemannian scalar curvatures of $g$ and $h$ are related the same kind of formula, where $b(n)$ should be 
replaced by $\frac{n-2}{4 (n-1)}$ and $\LL_{\theta}$ by the conformal Laplacian (and $\Delta_{\theta}$ by the Laplace-Beltrami operator). 
This analogy turns out to be very efficient: it is the key idea behind Schoen's 
proof (see Section \ref{Schoen}).

\subsubsection{The flat models}\label{flatmodels}

The \defini{standard CR sphere} $\ssp^{2n+1}$ (we will often omit the superscript)
is the unit sphere on $\mC^{n+1}$:
\[\ssp=\ensemble{(z_0,\dots,z_n)\in\mC^{n+1}}{\sum \abs{z_k}^2=1}\]
endowed with the corresponding CR structure. It is a strictly pseudoconvex
CR manifold; its automorphism group is $\Aut(S)=\PU{1,n+1}$, a finite 
quotient of $\SU{1,n+1}$. It is noncompact, connected and acts transitively on $\ssp$.

The \defini{Heisenberg group} is the CR noncompact manifold $\Heis$ obtained by 
removing one point of $\ssp$. It is therefore diffeomorphic to the Euclidean space
$\mR^{2n+1}$. Its automorphism group
is the stabilizer of the removed point in $\Aut(S)$, it acts non properly
and transitively and is connected.

These two CR manifolds are homogeneous and obviously locally isomorphic; they
are referred to as the flat models. They play the role of the Euclidean space in 
Riemannian geometry, or of the sphere and Euclidean space in conformal geometry.

For instance, there are local normal coordinates in any Riemannian manifold, where the metric is very close a Euclidean one.
There is an analogous local model for calibrated strictly-pseudoconvex manifolds: \cite{JL2} provides local ``normal'' coordinates 
in which the geometry is close to that of the Heisenberg group $\Heis$. In the Riemannian case, the local model (i.e. the Euclidean space) 
is global for simply connected complete flat manifolds. The following statement is the CR analogue.

\begin{prop}
A simply connected complete calibrated strictly-pseudoconvex CR manifold with 
vanishing curvature and torsion is CR equivalent to the
Heisenberg group. 
\end{prop}

Let us precise what ``complete'' means. The form $\theta$ being contact implies
that any two points in $M$ can be connected by a curve that is everywhere
tangent to the contact distribution. By minimizing the length of such curves,
one defines the \defini{Carnot distance} $d_\theta$. It is a genuine distance,
but does not derive from a Riemannian metric. By ``balls'' of $M$, we mean
balls with respect to the Carnot distance. A strictly pseudoconvex
CR manifold is said to be \defini{complete} if closed balls are compact.

\begin{defi}
We say that an open subset $U$ of a strictly pseudoconvex CR manifold $M$
is \defini{flat} if any $x$ in $U$ has a neighborhood which is CR isomorphic to an open subset
of $\ssp$.
\end{defi}

\subsection{Finite order rigidity}\label{forigidity}

For a general reference on local rigidity, see \cite{Kobayashi}, Theorems 3.2 
and 5.1. 

Let us start with the well-known rigidity of Riemannian geometry.
\begin{prop}\label{Riemannian}
A Riemannian metric on a manifold $M$ is rigid to order $1$, that is: two
isometries that have the same value and differential at some point are the same.
\end{prop}

In fact this result follows from a stronger statement. Let $OM$ be the
bundle of orthonormal frames on $M$ and $\Isom(M)$ its isometry group.
We look at its action on $OM$.
For each element ${\mathcal F}$ of the total space $OM$ (${\mathcal F}$
is thus the data of a point $x\in M$ and an orthonormal frame of $T_xM$), 
one defines the map
\begin{eqnarray}
  \Isom(M) &\to& OM \nonumber\\
  f &\longmapsto& f({\mathcal F}). \nonumber
\end{eqnarray}
Proposition \ref{Riemannian} asserts that this map is injective. In fact,
it is an embedding and its image is a closed submanifold of $OM$. The group
$\Isom(M)$, endowed with the corresponding differential structure, is a Lie group.

As a consequence, since the fibers of $OM$ are compact, the isometry group 
of a compact Riemannian manifold is compact. One even gets:

\begin{coro}\label{RiemannianC}
Let $U$ be an open set on a manifold $M$, $K\subset U$ be a compact set
with nonempty interior, $g$ be a Riemannian metric defined on $U$
and $G$ be a Lie group acting on $M$ and preserving $K$ and $g$. Then $G$
is compact.
\end{coro}

Now we turn to the rigidity of strictly pseudoconvex CR geometry.

\begin{prop}\label{rigidity}
Let $M$ be a strictly pseudoconvex CR manifold. The group $\Aut(M)$ is a 
Lie group and is rigid to order $2$, that is: if two automorphisms $f$, $f'$ 
have the same $2$-jet (the data of their derivatives up to order $2$) at some point, then
$f=f'$.
\end{prop}

As before, there is a principal bundle on $M$
in which $\Aut(M)$ embeds, but the fibers are no longer compact and
$\Aut(M)$ can thus be noncompact even when $M$ is compact.

As a direct consequence of Proposition \ref{rigidity}, two CR automorphisms of $M$
that coincide on an open set are the same.

The strict pseudoconvexity condition is of primary importance. For example, the product
$S^1\times\Sigma$ of the circle and any Riemann surface is a Leviflat CR manifold,
and the action of the infinite-dimensional diffeomorphism group of $S^1$ preserves
the CR structure.

\subsection{North-south dynamics}

The following result is a common feature of all ``rank $1$ parabolic geometries'',
that is of boundaries of negatively curved symmetric spaces. The standard CR
sphere $\ssp^{2n+1}$ is one of them: it bounds the complex hyperbolic space, seen as
the unit ball of $\mC^{n+1}$. Note that by an \defini{unbounded sequence} in
a topological space, we mean a sequence that is not contained in any compact set.

\begin{prop}\label{north-south}
Let $(\phi_k)_k$ be an unbounded sequence in $\Aut(\ssp)$. There exists a 
subsequence,
still denoted by $(\phi_k)_k$, and two
points (that may be the same) $p_+$ and $p_-$ on $\ssp$ such that:
\begin{eqnarray}
  \lim \phi_k(p) = p_+      &\quad&\forall p\neq p_- \\
  \lim \phi_k^{-1}(p) = p_- &\quad&\forall p\neq p_+
\end{eqnarray}
and the convergences are uniform on compact subsets of $\ssp\priv\{p_-\}$,
$\ssp\priv\{p_+\}$ respectively.

Moreover if the $\phi_k$'s are powers of a single automorphism $\phi$, then
$p_\pm$ are fixed point of $\phi$. The same result holds for a
noncompact flow, which has thus either one or two fixed points.
\end{prop}

An unbounded flow or automorphism\footnote{an automorphism is
said to be \defini{unbounded} if the sequence of its powers is unbounded}
is said to be \defini{parabolic} if it has
one fixed point, \defini{hyperbolic} if it has two of them. A bounded
flow or automorphism is said to be \defini{elliptic}.

\proof
The principle is to look at the action of $\Aut(\ssp)$ not only on the sphere
$\ssp$, but also in the complex hyperbolic space it bounds and on the projective
space $\CP^{n+1}$ it is embedded in.

\begin{figure}\begin{center}
  \input{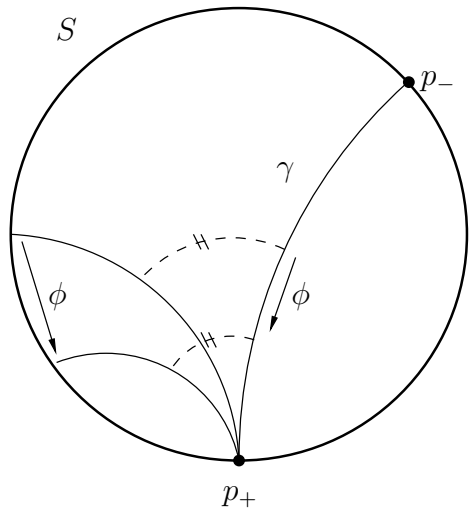}
  \caption{North-south dynamics.}\label{nsfigure}
\end{center}\end{figure}

The case when the $\phi_k$'s are powers of an automorphism $\phi$, or the
case of a flow, are simple linear algebra results. They are roughly described
by Figure \ref{nsfigure}, which shows the link between negative curvature of 
the hyperbolic spaces
and north-south dynamics: when a geodesic $\gamma$ is translated, any other 
geodesic is shrinked toward one of the ends of $\gamma$.

The general case can be deduced from the $KAK$ decomposition: every element
$\phi$ of the group $\Aut(\ssp)$ writes down as a product $\phi=k_1 a k_2$
where $k_1$, $k_2$ are elements of a maximal compact subgroup
$K\subset\Aut(\ssp)$ and $a$ is an element of a maximal noncompact closed abelian subgroup $A$.
The dimension of $A$ is the real rank of $\Aut(S)$, namely $1$. More precisely,
$A$ corresponds to a hyperbolic flow, that is a noncompact flow with two fixed
point on $S$ (one attractive, one repulsive). The general result
then follows from the compactness of $K$.
\finpreuve

%

\subsection{$(G,X)$-structures}\label{GX}

The notion of $(G,X)$-structure is a formalisation of Klein's geometry. In our
setting, they arise as a description of flat CR structures in term
of the model sphere $\ssp$ and its automorphism group.

Let $X$ be a manifold and $G$ a Lie group acting transitively on $X$. Assume 
that the action is \defini{analytic} in the following sense: an element that 
acts trivially on an open subset of $X$ acts trivially on the whole of $X$.
A \defini{$(G,X)$-structure} on a manifold $M$ is an atlas whose charts take 
their values in $X$ and whose changes of coordinates are restrictions of elements of
$G$. A diffeomorphism of $M$ is an \defini{automorphism} of its $(G,X)$-structure
if it reads in charts as restrictions of elements of $G$.

Let us consider the case when $G=\PU{1,n+1}$ and $X=\ssp$.
A flat strictly pseudoconvex manifold $M$ carries a 
$(G,X)$-structure  and its CR automorphisms
coincide with its $(G,X)$ automorphisms. Indeed, the flatness of $M$ 
means that it is locally equivalent to $\ssp$, thus it is sufficient to prove 
that any local automorphism of $\ssp$ can be extended into a global automorphism. 
This, in turn, follows from the order $2$ rigidity and the following fact:
any $2$-jet of a local automorphism of $\ssp$ can be realized as the $2$-jet of a 
global automorphism (see e.g. \cite{Spiro}).

The main tool we will need to study $(G,X)$-structures is the so-called
\defini{developping map}.
Let $M$ be a manifold endowed with a $(G,X)$-structure and $\tilde{M}$ be its
universal covering. Then there exists a differentiable map  
\[{\cal D}:\tilde{M}\to X\]
that is a local diffeomorphism and such that
for all automorphism $f$ of $\tilde{M}$, there exists some $\phi\in G$
satisfying
\begin{equation}
{\cal D}\rond f = \phi\rond {\cal D}.
\end{equation}
Note that in general, this developping map need not be a diffeomorphism 
onto his image, nor a covering map. It is unique, up to composition with 
an element of $G$. More details on $(G,X)$-structures can for example
be found in the classical \cite{Thurston}.

\section{Webster: a local Theorem}\label{Webster}

In 1977, Sydney Webster published the first work toward the Schoen--Webster
Theorem, \cite{Webster1}. 
Until the end of the paper, $M$ denotes a strictly pseudoconvex 
CR manifold of dimension $2n+1$.
\begin{theo}\label{Wlocal}
If $M$ is compact and $\Aut_0(M)$ is noncompact, then $M$ is flat.
\end{theo}

There are several reasons why this result has raised a lot of efforts to be improved.
First, it is a local statement though Webster gave in the same paper
a very specific global result:
\begin{theo}\label{Wglobal}
If $M$ is compact and has finite fundamental group and $\Aut_0(M)$
is noncompact, then $M$ is globally equivalent to the standard sphere.
\end{theo}

Second, he assumes that $M$ is compact and that the 
\emph{identity component} $\Aut_0(M)$ is noncompact. We refer
to these hypotheses as the \defini{compactness assumption} and the \defini{connectedness 
assumption}.

His paper also contains a result on connected groups of CR automorphisms 
having a fixed point we shall discuss briefly. 
\begin{theo}\label{Wfixed}
If $M$ is compact and $\Aut_0(M)$ admits a noncompact one-parameter Lie subgroup 
$G_1$ that has a fixed point $p_0$, then $M$ is globally equivalent to
the standard sphere $\ssp$.
\end{theo}

Let us turn to the proofs of these three results.

\subsection{Canonical calibration}

The following result is the key to the local statement.

\begin{lemm}\label{relative_invariant}
For each calibration $\theta$ on $M$ there is a continuous nonnegative 
function $F_\theta$ on $M$ such that:
\begin{enumerate}
\item $F_\theta$ vanishes on a given open set $U$ if and only if $U$ is flat,
\item on the open set where $F_\theta$ is positive, it is smooth,
\item the family $(F_\theta)_\theta$ is homogeneous of degre $-1$ : 
      \begin{equation}
        F_{\lambda\theta}=\abs{\lambda}^{-1}F_\theta.
      \end{equation}
\end{enumerate}
\end{lemm}
Such a family of functions $(F_\theta)_\theta$ is called a \defini{relative
invariant} after Cartan's one (see the proof below). 
Most of the time, the are given by the norm of a curvature
tensor.

A point where $F_\theta$ vanishes for some (thus for all) calibration $\theta$
is said to be \defini{umbilic}.

Let us show the interest of such functions.
Pick any calibration $\theta$ of $M$ whose Levi form is positive and define
\begin{equation}
\theta^*=F_\theta \theta.
\end{equation}
Then $\theta^*$ is a continuous 1-form that vanishes on the flat part
of $M$ and is a smooth calibration everywhere else. It is canonical,
for if $\theta'=\lambda\theta$ is another calibration with $\lambda>0$
(that is, whose Levi form is positive),
\begin{eqnarray}
\theta'^* &=& F_{\lambda\theta} \lambda\theta \nonumber\\
          &=& \lambda^{-1}F_\theta \lambda\theta \nonumber \\
          &=& \theta^* \nonumber.
\end{eqnarray}
We call $\theta^*$ the \defini{canonical calibration} of $M$ although it is not
a genuine calibration unless $M$ contains no umbilic points. If $M$ is flat,
$\theta^*$ is zero and, therefore, useless.

\proof
Lemma \ref{relative_invariant} follows from the study of invariants of
calibrated CR manifolds.

If $n>1$, one can derive from the Chern-Mother
curvature a tensor $S$ on some bundle $T$ over $M$ that only depends upon 
the CR structure and vanishes on an open set $U$ if and only if $U$ is 
flat. A calibration $\theta$ induces, \lat{via} the Levi form,
a metric on $T$. The corresponding norm $\norme{S}_\theta$ 
of $S$ yields the desired function. See \cite[page 201]{Burns-Schnider}
or \cite[page 35]{Webster2} for details.

If $n=1$, $S$ is always zero even when $M$ is not flat so that
Cartan's relative invariant is needed. It is a function
$r_\theta$ on $M$, associated to a calibration $\theta$, that vanishes on an
open set if and only if it is flat; the family $(r_\theta)$ is homogeneous of 
order $-2$,
thus $F_\theta=\sqrt{r}$ does the job. For details, one can look at
\'{E}lie Cartan's work \cite{Cartan1,Cartan2} or, for a more modern presentation,
at the book of Howard Jacobowitz \cite{Jacobowitz}.
\finpreuve

\subsection{The local theorem}

Let us give an outline of the proof of Theorem \ref{Wlocal} given by 
Webster. We shall see later that a stronger statement 
can be proved with the same tools.

\proof
Assume $M$ is not flat; we will show that any one-parameter subgroup
of $\Aut_0(M)$ has a compact closure, which implies the compactness
of $\Aut_0(M)$ by a theorem of Deane Montgomery and Leo Zippin
\cite{Montgomery-Zippin}.

Let $G_1$ be a nontrivial one-parameter subgroup of $\Aut_0(M)$ with
infinitesimal generator $Y$ on $M$. Choose some calibration $\theta$; by assumption
$F_\theta$ is positive on an open set $U$. Since the vanishing of
$F_\theta$ is independent of $\theta$, $U$ is invariant under the flow of $Y$.

Consider the function $\eta=\theta^*(Y)$, on $U$. Assume it vanishes identically.
Then $Y$ lies in the contact distribution. Moreover, since $\theta^*$ is a CR 
invariant form, ${\mathcal L}_Y\theta^*$ vanishes. Cartan's magic formula yields:
\[0 = {\mathcal L}_Y\theta^*=Y\contraction d\theta^*+d\eta = Y\contraction d\theta^*,\]
so $Y$ is identically zero, which contradicts the order two rigidity.

We may therefore assume that $\eta>0$ somewhere, replacing $Y$ by $-Y$ if
necessary. Choosing $\varepsilon$ sufficiently small, the set
$U_\varepsilon$ defined by the inequation $\eta(p)\geqslant \varepsilon$
has non empty interior. It is closed in $M$, thus is compact, and is
invariant under the flow of $Y$. 

The closure $\adherence{G_1}$ of $G_1$ in $\Aut_0(M)$ is a Lie group that
preserves the compact $U_\varepsilon$ and the Webster metric of $\theta^*$ on it,
thus is compact (Corollary \ref{RiemannianC}).
\finpreuve

\subsection{The global result}

Webster derives Theorem \ref{Wglobal} from a weak form
of Proposition \ref{north-south} and a (now) standard use
of $(G,X)$-structures.

\proof By Theorem \ref{Wlocal}, $M$ and its universal covering $\tilde{M}$ are
flat. Therefore they can be developped as $(\SU{1,n+1},\ssp)$-structures. Since $M$ 
has finite fundamental group, $\tilde{M}$ is compact and the developping
map ${\mathcal D}:\tilde{M}\rightarrow \ssp$ is a covering map. But $\ssp$
admits no nontrivial covering and $\tilde{M}$ is globally equivalent to $\ssp$.
By Montgomery-Zippin Theorem \cite{Montgomery-Zippin}, there exists
some closed noncompact one-parameter subgroup $G_1$ of $\Aut_0(M)$. This group
lifts to a one parameter subgroup $\tilde{G_1}$ acting on $\tilde{M}=\ssp$.

>From Proposition \ref{north-south}, we know that $\tilde{G_1}$
has either one or two fixed points. In both cases $G_1$ has at least a fixed point.

Let $p$ be a fixed point of $G_1$. The lifts of $p$ are fixed points of $\tilde{G_1}$
of the same type (attractive, repulsive or both). But $\tilde{G_1}$ has at most
one fixed point of a given type, thus $\tilde{M}\rightarrow M$ must be a one-sheeted
covering.
\finpreuve

\subsection{One-parameter subgroups with a fixed point}

Theorem \ref{Wfixed} is based on a principle of extension of local conjugacy,
making use of the dynamics on the model space. We detail a similar argument at the
end of the paper, using it in the proof of the compact case of the Schoen--Webster
Theorem.

\proof
Let $Y$ be an infinitesimal generator of $G_1$. According to Theorem \ref{Wlocal}, $M$ is flat 
so there is an isomorphism between a neighborhood $U$ of $p_0$ and an open set $U'$ of $\ssp$.
Denote by $Y'$ the vector field on $U'$ corresponding to the restriction of $Y$ to $U$.
Then $Y'$ extends uniquely to a CR vector field on $S$, which has a fixed point $p'_0$.

If $Y'$ is elliptic, then it follows from the finite order rigidity that $G_1$ is compact,
in contradiction with the assumptions.
If $Y'$ is parabolic, then one can use it to extend the conjugacy between $U$ and $U'$ to the basins of 
attraction and repulsion of $p'_0$, therefore $M$ is globally equivalent to $\ssp$.
If $Y'$ is hyperbolic, the same argument shows that there is an open set $V\subset M$
that is conjugate to either $\ssp$ or $\ssp$ with a point (namely the second fixed point of $Y'$) deleted.
Since $Y$ is a complete vector field with isolated zeros, $M$ itself must be globally 
equivalent to either $\ssp$ or $\ssp$ with a point deleted.
\finpreuve

\section{Kamishima and Lee: two ways from local flatness to global rigidity}

Yoshinobu Kamishima seems to be the first to prove the local to global
statement (under both the compactness and connectedness assumptions)
in a workshop in honor of Obata held at Keio University in 1991.
He announced the result in the proceedings \cite{Kamishima1} and the
complete proof appeared a while after \cite{Kamishima2}.

\begin{theo}\label{Kglobal}
If $M$ is flat, compact and $\Aut_0(M)$ is noncompact then
$M$ is globally equivalent to the standard sphere.
\end{theo}

We will not detail his proof at all, but let us quote an interesting corollary
he gave in relation with the so-called \defini{Seifert conjecture}. This celebrated
conjecture states that any non singular vector field on the $3$-dimensional
sphere has at least one closed orbit. It was disproved
for $\diffb{1}$ vector fields by Paul Schweitzer \cite{Schweitzer}, then in
$\diffb{\infty}$ regularity by Krystyna Kuperberg \cite{Kuperberg}. 
The question was then raised for vector fields preserving some geometric
structure. Kamishima's following result gives an answer for vector fields
preserving a CR structure.

\begin{coro}
If $M$ is a rational homology sphere endowed with a strictly pseudoconvex
CR structure, then any nonsingular CR vector field
on $M$ has a closed orbit.
\end{coro}

In \cite{Lee}, John Lee proved Theorem \ref{Kglobal} independently of 
Ka\-mi\-shi\-ma.
His method relies on Webster's Theorem \ref{Wfixed}: he proves
\begin{theo}\label{Lglobal}
If $M$ is compact and $\Aut_0(M)$ admit a closed noncompact one-parameter 
subgroup $G_1$, then $G_1$ has a fixed point.
\end{theo}

Once again, the Montgomery--Zippin Theorem is used to deduce Theorem \ref{Kglobal}
from Theorems \ref{Lglobal}, \ref{Wlocal} and \ref{Wfixed}.

\proof
Let $Y$ be an infinitesimal generator of $G_1$ and assume by contradiction
that $Y$ has no zero on $M$.

Note first that $Y$ must be somewhere tangent to the contact distribution $\xi$:
otherwise $Y$ would be the Reeb vector field of a unique calibration, thus would
preserve the associated Webster metric. 

The first part of the proof consists in 
understanding the set of points where $Y\in\xi$; it is a classical computation
that involves only the contact structure on $M$:
pick any calibration $\theta$ of $M$ and define $\eta=\theta(Y)$.
Then one can show, using that $Y$ has no zero, that $0$ is a regular value of $\eta$.
Therefore $H=\{\eta=0\}\subset M$
is a nonempty, compact, embedded hypersurface along which $Y$ is tangent 
to both $H$ and $\xi$.

The next step consists in proving that one can find a new calibration
such that ${\mathcal L}_Y \theta=0$ and ${\mathcal L}_Y d\theta=0$
at every point of $H$. It is easy to get the first condition by rescaling
$\theta$; then a rather tedious computation allows Lee to
refine the rescaling in order to get the second condition.

These two conditions imply that the Webster metric on $TM$ is preserved
by the flow of $Y$ along $H$. It follows for any sequence $(f_i)$ in
$G_1$, $(f_i\restric{H})$ converges in $\diffb{\infty}$ topology.
Using the complex operator $J$, it is then possible to prove
that the $2$-jets of the sequence $(f_i)$ converge at all points of $H$.
By the order $2$ rigidity, $(f_i)$ is convergent in $G_1$, a contradiction.
\finpreuve

\section{Schoen: Yamabe problem methods}\label{Schoen}

The aim of this section is to survey the proof of Schoen--Webster theorem by R. Schoen in \cite{Schoen}. For convenience, we only deal 
with the compact case, even though \cite{Schoen} also considers the non-compact case with the same kind of techniques, based on global 
analysis. R. Schoen first proves that the conformal group is compact for any closed Riemannian manifold which is not conformally equivalent 
to the standard sphere. Then he explains how to adapt the proof in a CR setting, which is what we want to develop below. Another proof 
of the conformal group compactness is given in \cite{Heb} : it is a bit shorter but relies on the positive mass theorem, which makes it less
elementary than what follows.  

\subsection{Yamabe theorem.} The celebrated Yamabe problem is basic in conformal geometry: is there a metric with constant scalar curvature in each conformal class 
of a given closed manifold ? This question was the beginning of a long story: see the excellent \cite{Heb} or \cite{LP} for an exhaustive 
account. The answer to the problem is yes and the proof relies on a careful study of the conformal Laplacian. 

As explained in \cite{JL1}, there is a deep analogy between conformal and CR geometry. In particular, Yamabe theory has a counterpart in 
the CR realm, which enables R. Schoen to extend his conformal geometry arguments to the CR case. 

In order to develop a Yamabe theory in the CR setting, one needs a Sobolev-like analysis. In the conformal case, the natural conformal 
operator is elliptic, so that its analysis is rather standard. In the CR case, the corresponding natural operator $\LL_\theta$ is only 
subelliptic. G. Folland and E. Stein \cite{FS} (see also paragraph 5 of \cite{JL1}) have nonetheless developped a powerful theory which 
yields the necessary tools. 

As in conformal geometry, given a calibration $\theta$, we define the CR Yamabe invariant $Q(M,\theta)$ as the infimum of the functional 
$$
\int_M \phi \LL_\theta \phi \; \theta \wedge d\theta^n
$$
over the elements $\phi$ of the unit sphere in the Lebesgue space $L^{\frac{2n+2}{n}}(M)$. The choice of this exponent is related to the
transformation law for the volume form: if $\theta' = u^{\frac{2}{n}} \theta$ for some positive function $u$, then
$$
\theta' \wedge (d\theta')^n = u^{\frac{2n+2}{n}} \theta \wedge d\theta^n.
$$
It turns out that $Q(M,\theta)$ is a CR invariant. 

D. Jerison, J. M. Lee \cite{JL1}, N. Gamara and R. Yacoub \cite{Gam}, \cite{GY} adapted the proof of the conformal Yamabe 
theorem to prove the 
\begin{theo}\label{yamabeCR}
A closed strictly pseudoconvex CR manifold admits a calibration with constant 
scalar curvature $1$ (resp. $0$ and 
$-1$) if its CR Yamabe invariant is positive (resp. zero and negative). 
\end{theo}

We will only need the nonpositive (and easiest) case, which was settled by \cite{JL1}. 

\subsection{The proof.} Theorem \ref{yamabeCR} leads to the

\begin{prop}\label{yamabenegatif}
When the CR Yamabe invariant is nonpositive, the CR automorphism group is compact.  
\end{prop}

\proof
We prove that CR automorphisms are isometries for the Webster metric of a calibration; since the isometry group of a closed 
Riemannian manifold is compact, the result will follow. Endow $M$ with a calibration $\theta$.

If $Q(M)=0$, we can assume $\theta$ has vanishing scalar curvature (Yamabe). A CR automorphism $F$ of $M$ then obeys 
$F^*\theta = u^{\frac{2}{n}} \theta$ with $\LL_\theta u = \Delta_\theta u = 0$ ($F^*\theta$ has scalar curvature $F^*R_\theta = 0$). An 
integration by parts yields 
$$
0=\int u \Delta_\theta u = \int L_\theta \left(du\restric{\xi},du\restric{\xi} \right) \theta \wedge d\theta^n.
$$ 
So we can write $du = f \theta$, which implies $0=df \wedge \theta + f d\theta$. Since $d\theta$ is definite on the kernel 
$\xi$ of $\theta$, $f$ vanishes, so $u$ is constant. Since 
$$
\vol_\theta(M)= \vol_\theta(F(M)) = \vol_{F^*\theta}(M)= u^{\frac{2n+2}{n}} \vol_\theta(M),
$$ 
$u$ is constant to $1$: $F$ preserves $\theta$ hence $W_\theta$. 

If $Q(M)<0$, we can make a similar argument: we are left to show that a solution $u$ of  
$$
\Delta_\theta u = b(n) \left( u - u^{\frac{n+2}{n}} \right) 
$$ 
is constant to $1$. It follows from a weak maximum principle. At a maximum point, $\Delta_\theta u$ is nonnegative so that the equation 
ensures $u\leq 1$. At a minimum point, one finds $u\geq 1$ for the same reason. Therefore $u$ is constant to $1$. 
\finpreuve

The following lemma is the key to complete the proof. We denote by $D_r$
the ball of radius $r$ in $\mR^{2n+1}$. To avoid technical details, we do not give the precise statement (cf. \cite{Schoen}).  

\begin{lemm}\label{bouleplate}
Let $F:(D_1,\theta)\to(N,\sigma)$ be a CR diffeomorphism. We assume $\theta$ is close to the Heisenberg calibration and $\sigma$ 
has vanishing scalar curvature. If $\lambda := \sqrt{(F^*\sigma/\theta)(0)}$ denotes the dilation factor at $0$, then:  
\begin{itemize}
\item the dilation factor is almost constant to $\lambda$, i.e. $F^*\sigma/\theta \approx \lambda$ on $D_{1/2}$ ; 
\item images of balls have moderate eccentricity, i.e. $F(D_{1/2}) \approx B( F(0) , \lambda/2 )$ ; 
\item the total curvature and the torsion are small when the dilation factor is large, i.e. $\abs{\Rm_\sigma} \lesssim \lambda^{-2}$ and 
$\abs{\T_\sigma} \lesssim \lambda^{-2}$ on $B(F(0) , \lambda/2 )$.
\end{itemize}
\end{lemm}
 
\proof
By scaling $\sigma$, we can assume $\lambda=1$. Write $F^*\sigma=u^{\frac{2}{n}} \theta $ and observe that $\R_\sigma = 0$ implies 
$\LL_\theta u=0$. Since $\theta$ is close to the Heisenberg calibration, $\LL_\theta$ is close to the Heisenberg subelliptic Laplacian, so 
that $u$ satisfies a Harnack inequality (\cite{JL1}, 5.12) : $\sup u \leq C \inf u$, with a controlled constant. The first and second 
assertions follow. Subelliptic regularity (\cite{JL1}, 5.7) also yields a $C^2$ bound on $u$, hence the third assertion.      
\finpreuve

Now we can finish the proof of the

\begin{theo}[Schoen--Webster]
The CR automorphism group of a closed strictly pseudo-convex CR manifold which is not CR equivalent to a standard sphere is compact.
\end{theo}

Here, we only deal with $C^0$ topology. Thanks to a bootstrap argument, \cite{Schoen} proves that all $C^k$ topologies, $k \geq 0$, are the 
same. They also coincide with the Lie group topology. 
 
\proof
Assume $M^{2n+1}$ is a closed strictly pseudo-convex CR manifold with non-compact conformal group and choose a calibration $\theta$. Ascoli theorem yields CR 
automorphisms $F_i$ and points $x_i$ such that the dilation factors 
$$
\lambda_i := \sqrt{(F_i^*\theta/\theta)(x_i)} = \max \sqrt{(F_i^*\theta/\theta)}
$$ 
go to infinity. 

The rough idea of the proof consists in multiplying the calibration $\theta$ by suitable Green functions so as to build a sequence of
conformal scalar flat blow ups; then lemma \ref{bouleplate} will enable us to find a sequence of larger and larger balls endowed with a 
calibration of smaller and smaller curvature and torsion: taking a limit, we will realize $M$ minus a point as a Heisenberg group; a 
last effort will seal the fate of the missing point.

To begin with, we can choose a small $\epsilon>0$ such that the geometry of all the balls of radius $\epsilon$ in $M$ is close to that 
of the Heisenberg group. Then we choose points $y_i$ outside $F_i(B(x_i,\epsilon))$ and use a standard trick in Yamabe theory. Since the 
CR Yamabe invariant is positive (\ref{yamabenegatif}), the operator $\LL_\theta$ is positive. Therefore, there are Green functions $G_i$,
i.e. preimages of Dirac distributions $\delta_{y_i}$ (cf. \cite{Gam} for instance): outside $y_i$, they are smooth, satisfy $\LL_\theta G_i=0$ 
and we can normalize them so that their minimum value is $1$. Put $z_i := F_i(x_i)$ and consider the calibration 
$$
\theta_i := \left( \frac{G_i}{G_i(z_i)} \right)^{\frac{2}{n}} \theta,
$$
defined outside $y_i$. It has vanishing scalar curvature. 

We can assume $y_i$ converges to $y$, $z_i$ converges to $z$ and $G_i$ converges to $G$ on compact sets of $M \priv \{y\}$. Besides, one can
show that $G_i(z_i)$ remains bounded, that is $y \not= z$: it stems from a convenient use of lemma \ref{bouleplate} and from a Harnack 
inequality for the dilation factor between $F_i^*\theta_i$ and $\theta$. So we can assume $G_i(z_i)$ converges. 


Therefore $\theta_i$ tends to a calibration $\theta_\infty = c G^{\frac{2}{n}} \theta$ on the compact sets of $M \priv \{y\}$. Now 
lemma \ref{bouleplate} ensures that, roughly, $\theta_i$ has curvature and torsion of magnitude $\lambda_i^{-2}$ on 
$F_i(B_\theta(x_i,\epsilon/2)) \approx B_{\theta_i}(z_i,\lambda_i \epsilon/2)$, so that letting $i$ go to infinity, we conclude our manifold, 
outside $y$, is CR equivalent to a calibrated strictly pseudo-convex CR manifold with vanishing curvature and torsion; and it happens to be 
complete and simply connected (it is a nondecreasing union of topological balls), so that it is $\Heis^{2n+1}$.  
 
Thus there is a CR diffeomorphism $F$ between $M$ minus $y$ and the standard sphere minus some point, $\infty$. In the neighbourhood of 
$\infty$ in $\ssp^{2n+1}$, consider a CR equivalent Heisenberg calibration $\sigma$. Writing $F^*\sigma=u^{\frac{2}{n}} \theta$, we obtain 
$\LL_\theta u=0$ outside $y$, since $\sigma$ has vanishing scalar curvature. Extending $F$ at $y$ amounts to show that $u$ has a removable 
singularity at $y$. But the integral of $u^{\frac{2n+2}{n}}$ over some ball is exactly the volume of the image of this ball through $F$, 
which is bounded by the volume of the standard sphere; it follows (proposition 5.17 in \cite{JL1}) that $u$ is a weak solution of the 
equation $\LL_\theta u=0$ over a neighborhood of $y$ so that it extends as a smooth function in the neighborhood of $y$ (5.10, 
5.15 in \cite{JL1}). Thus $(M,\theta)$ is CR equivalent to the standard sphere. 
\finpreuve

\section{Frances: a unified dynamical proof}

Charles Frances recently gave a unified proof of the Ferrand-Obata
and Schoen--Webster Theorems \cite{Charlie2}; in fact he also proves analoguous 
results for quaternionic-contact and octonionic-contact geometries.

To obtain these results, he uses the setting of Cartan geometries (see \cite{Sharpe} for
a detailed account on this topic). Given a model homogeneous space $X=G/P$, a 
\defini{Cartan geometry modelled on} $X$ on a manifold $M$ consists of:
\begin{itemize}
\item a $P$-principal bundle $B\rightarrow M$ and
\item a $1$-form $\omega$ on the total space $B$ with values in the Lie
      algebra $\al{g}$.
\end{itemize}
The form $\omega$ is called the \defini{Cartan connection} of the structure and is
supposed to satisfy some compatibility conditions we do not detail.

The points of $B$ play the role of  ``adapted'' frames (like orthonormal frames
for Riemannian geometry). The Cartan connection is used to identify
infinitesimally $B$ with $G$: in particular, it is asked that at each point
$p\in B$, $\omega_p$ is an isomorphism between $T_pB$ and $\al{g}$.

The geometries Frances is concerned with are modelled on the homogeneous
spaces $\partial\KH^d=G/P$ where $\KH^d$ is the hyperbolic space based on
$\mK=\mR$, $\mC$, $\mH$ or $\mO$, $G$ is the isometry group of
$\KH$ and $P$ is the stabilizer of a boundary point. Note that when
$\mK=\mC$, $X=\ssp$.

For each of these Cartan geometries, the ``equivalence problem'' has been
solved, that is: there exists a construction that gives for any conformal,
strictly pseudoconvex CR, etc. structure on $M$ a corresponding Cartan
structure $B,\omega$ such that isomorphisms of the original structure induce
isomorphisms of the Cartan structure and reciprocally. The Cartan
structure is not unique, one can impose further assumptions. In particular
the Cartan connection can be chosen ``regular'' (a technical condition
involving the curvature of $\omega$) for the geometries considered here.

We can now state the result of Frances.
\begin{theo}
Let $(M,B,\omega)$ be a Cartan geometry modelled on $X=\partial\KH^d$,
with regular connection. If $\Aut(M,\omega)$ acts nonproperly on $M$,
then $M$ is isomorphic to either $X$ or $X$ with
a point deleted.
\end{theo}

\proof
The first and main step is to prove that any sequence $(f_k)$ of automorphisms
of $M$ that acts nonproperly admit a subsequence that shrinks an open set
$U\subset M$ onto a point $p$. The principle is to use some sort of 
developping map from the space of curves
on $M$ passing through $p$ to the space of curves on $X$ passing through
a given base point $o$. Then, choosing an appropriate family of curves
in the model and its north-south dynamics, one gets the desired property 
on $M$.

Then one proves that an open set that collapses to a point must be
flat. Note that in the CR case, one could use the Webster metric of the canonical
calibration (see Section \ref{Webster}).
As a consequence, one can choose $U$ to be of the form
\[U=\Gamma\backslash (X\priv\{o\})\]
where $\Gamma$ is a discrete subgroup of the stabilizer $P$ of $o\in X$.

The final step is a result on geometrical rigidity of embeddings:
if a flat manifold
$\Gamma\backslash (X\priv\{o\})$ embeds in $M$, then either
$M=\Gamma\backslash (X\priv\{o\})$ or $\Gamma=\{\id\}$. In the latter case,
$M=X$ or $M=X\priv\{o\}$.

The conclusion follows since the automorphism group
of $\Gamma\backslash(X\priv\{o\})$ acts properly when $\Gamma$
is not trivial.
\finpreuve

\section{Gathering a geometric proof in the compact case}\label{gathering}

In this last section, we give a geometric proof of the Schoen--Webster Theorem
under the compactness assumption.
It is not elementary, as it makes use of Lemma 
\ref{relative_invariant}.
However: it is a geometric proof, thus gives an alternative to Schoen's techniques;
it do not rely on the Montgomery-Zippin Theorem, holds without the connectedness 
assumption and is quite short, which makes it an improvement of those of
Webster, Ka\-mi\-shi\-ma and Lee together.

It does not pretend to originality, since it relies on arguments of
Webster \cite{Webster1} and Frances and Tarquini \cite{Charlie1}, rephrased.

\subsection{The local statement}

\begin{theo}
If $M$ is compact and $\Aut(M)$ is noncompact, then $M$ is flat.
\end{theo}

\proof
Suppose that $M$ is not flat; then the canonical calibration $\theta^*$
defined thanks to Lemma \ref{relative_invariant} does not vanish identically. 
Denote by $W$
the Webster metric associated with $\theta^*$: it is continuous on $M$,
smooth and positive definite on the open set $U$ of nonumbilic points
and zero on its complementary $F$. For all $x$ and $y$ in $M$ let
\[d(x,y)=\inf_\gamma \int_\gamma \sqrt{W(\dot{\gamma})} \]
define the natural semimetric associated to $W$ (not to be confused
with the Carnot metric of Section \ref{flatmodels} : here the infimum is
taken on all curves connecting $x$ to $y$). 
We have $d(x,y)=0$ if and only
if $x$ and $y$ are in $F$, in particular $d$ is a genuine metric on $U$.

If $U=M$, then $\Aut(M)$ preserves a Riemannian metric, thus is compact. Otherwise,
$F$ is nonempty, the distance $d(x,F)$ is finite for every $x\in M$ and 
we can define the set
$U_\varepsilon=\ensemble{x\in U}{d(x,F)\geqslant\varepsilon}$
for any positive $\varepsilon$. This set 
is compact and has nonempty interior for $\varepsilon$ small enough.

Now $\Aut(M)$ preserves $U_\varepsilon$ and its Webster metric, thus is
compact.
\finpreuve

\subsection{The local-to-global statement}

\begin{theo}\label{loctoglob}
If $M$ is flat and $\Aut(M)$ acts nonproperly,
then $M$ is globally equivalent to the standard CR sphere $\ssp$ or to
$\ssp$ with a point deleted.
\end{theo}

This result follows, by a principle of ``extension of local conjugacy'',
from the dynamics of unbounded sequences of $\Aut(\ssp)$. Note that we do not
use the compactness assumption for this part.

The end of the section is dedicated to the proof of Theorem \ref{loctoglob}.
Note that it holds as it is for any ``rank-one parabolic'' $(G,X)$-structure 
(namely $X=\partial\KH^n$ with $\mK=\mR$, $\mC$, $\mH$ or $\mO$).

\subsubsection{Set up: developping the dynamics.}

We assume that $M$ is flat, thus carries a $(\SU{1,n+1},\ssp)$-structure,
and that $\Aut(M)$ acts nonproperly: there is a convergent sequence $x_i\in M$
and a sequence $f_i\in\Aut(M)$ going to infinity (that is, having no convergent 
subsequence), such that $y_i=f_i(x_i)$ converges in $M$. We set $x_\infty=\lim x_i$ and
$y_\infty=\lim y_i$. 

Let $\tilde{M}$ be the universal cover of $M$. There are lifts 
$(\tilde{x}_i)_{i\in\mN\cup\{\infty\}}$, $(\tilde{y}_i)_{i\in\mN\cup\{\infty\}}$ 
and $\tilde{f}_i$ such that $\lim \tilde{x}_i = \tilde{x}_\infty$,
$\lim \tilde{y}_i = \tilde{y}_\infty$ and $\tilde{y}_i=\tilde{f}_i(\tilde{x_i})$.
Moreover, the sequence $(\tilde{f}_i)$ has no convergent subsequence in $\Aut(\tilde{M})$.

Let ${\cal D}:\tilde{M}\rightarrow \ssp$ be the developping map of $M$ and $\phi_i$
be a sequence of $\Aut(\ssp)$ such that ${\cal D}\tilde{f}_i=\phi_i{\cal D}$.
If $(\phi_i)$ had a convergent subsequence, by the order $2$ rigidity and since
$\phi_i$ and $\tilde{f}_i$ are locally conjugated, so would $(\tilde{f}_i)$. Thus
$(\phi_i)$ is unbounded and admit a North-South dynamics, whose poles are denoted
by $p_+$ and $p_-$.

Since ${\cal D}(\tilde{y}_i)=\phi_i{\cal D}(\tilde{x}_i)$, we have either
${\cal D}(\tilde{y}_\infty)=p_+$ or ${\cal D}(\tilde{x}_\infty)=p_-$. Up to inverting 
the $f_i$'s and exchanging
the $x_i$'s and the $y_i$'s, we assume that ${\cal D}(\tilde{y}_\infty)=p_+$.

\subsubsection{Stretching injectivity domains.}

A subset of $\tilde{M}$ is said to
be an \defini{injectivity domain} if the developping map is one-to-one on
its closure.

We denote by $U_0$ an open connected injectivity domain containing $\tilde{y}_\infty$
and we let $V_0={\cal D}(U_0)$. We choose an open connected injectivity domain $\Omega$ 
containing $\tilde{x}_\infty$ and having connected boundary $\bound\Omega$  
whose image ${\cal D}(\bound\Omega)$ does not contain $p_-$.
Up to extracting a subsequence, we can assume that for all $i$,
$\tilde{x}_i\in\Omega$ and $\tilde{y}_i\in U_0$.

According to Proposition \ref{north-south}, there is an increasing sequence of open sets
$V_i\subset \ssp$ ($i>0$) such that, extracting a subsequence if necessary:
\begin{enumerate}
\item for all $i$, ${\cal D}(\bound\Omega)\subset V_i$,
\item $\bigcup V_i = S\priv\{p_-\}$,
\item for all $i$, $\phi_i(V_i)\subset V_0$.
\end{enumerate}

Let $\delta:U_0\rightarrow V_0$ be the restriction of ${\cal D}$ and define the following
open connected injectivity domains: $U_i=\tilde{f}_i^{-1}\rond\delta^{-1}\rond\phi_i(V_i)$.
Since we assumed $\tilde{x}_i\in\Omega$ and $\tilde{y}_i\in U_0$, we get
\begin{equation}
U_i\cap\Omega\neq\varnothing\quad\forall i
\end{equation}
and by construction we have
\begin{equation}
{\cal D}(\bound \Omega)\subset{\cal D}(U_i)=V_i\subset{\cal D}(U_{i+1})=V_{i+1}\quad\forall i.
\end{equation}

\subsubsection{Monotony and consequences}

We prove that $(U_i)$ (or a subsequence) is an increasing sequence. 

If we can extract a subsequence such that $U_i\subset\Omega$ for all $i$,
since $\Omega$ is an injectivity domain and $({\cal D} U_i)$ is increasing,
$(U_i)$ must be increasing.

Otherwise we use the following Lemma.
\begin{lemm}\label{interinject}
Let $A$, $B$ be two injectivity domains such that $B$ is open,
$A$ is connected and $A\cap B\neq\varnothing$.
If ${\cal D}(A)\subset{\cal D}(B)$, then $A\subset B$.
\end{lemm}

\proof
Since $A$ is connected,
we only have to prove that $A\cap B$ is open and closed in $A$.

First, $B$ is open so that $A\cap B$ is open in $A$.
Second, let $y$ be a point in $A\cap \adherence{B}$. Since ${\cal D}(A)\subset{\cal D}(B)$,
there is a $z\in B$ such that ${\cal D}(z)={\cal D}(y)$. But since $B$ is an injectivity
domain and $y$ belongs to the closure of $B$, $z=y$ and $y\in B$. Therefore,
$A\cap B$ is closed in $A$.
\finpreuve

When $U_i\not\subset\Omega$, $\bound\Omega\cap U_i\neq\varnothing$ and we can apply 
Lemma \ref{interinject}:
$\bound\Omega\subset U_i$. But then $U_i\cap U_{i+1}\neq\varnothing$
thus by the same argument: $U_i\subset U_{i+1}$.

Now let $U_\infty=\bigcup U_n$; ${\cal D}$ is a diffeomorphism from $U_\infty$ to
$\ssp\priv\{p_-\}$.

If $U_\infty\neq\tilde{M}$, let $x$ be a point of the boundary
of $U_\infty$. If ${\cal D}(x)$ were in 
$\ssp\priv\{p_-\}$, for any neighborhood $W$ of $x$, 
${\cal D}(W\cap U_\infty)$ would meet any neighborhood
of any inverse image of  ${\cal D}(x)$, contradicting the injectivity of ${\cal D}$
on $U_\infty$. Therefore the boundary of $U_\infty$ consists of $x$ alone
and ${\cal D}$ is a global diffeomorphism from $\tilde{M}=U_\infty\cup\{x\}$ 
to $\ssp$.

We thus proved that $\tilde{M}$ is equivalent to either $\ssp$ or $\ssp\priv\{p_-\}$. Moreover,
any inverse image of $y_\infty$ in $\tilde{M}$ is an attracting point, thus is
$\tilde{y}_\infty$: $\tilde{M}$ is one-sheeted and
$M$ is itself equivalent to either $\ssp$ or $\ssp\priv\{p_-\}$.

\bibliographystyle{alpha}
\bibliography{biblio}

\begin{thebibliography}{Car32b}

\bibitem[BS76]{Burns-Schnider}
D.~Burns, Jr. and S.~Shnider.
\newblock Spherical hypersurfaces in complex manifolds.
\newblock {\em Invent. Math.}, 33(3):223--246, 1976.

\bibitem[Car32a]{Cartan1}
\'{E}lie Cartan.
\newblock Sur la g\'eom\'etrie pseudo-conforme des hypersurfaces de l'espace de
  deux variables complexes.
\newblock {\em Ann. Mat. Pura Appl. (4)}, 11:17--90, 1932.

\bibitem[Car32b]{Cartan2}
\'{E}lie Cartan.
\newblock Sur la g\'eom\'etrie pseudo-conforme des hypersurfaces de l'espace de
  deux variables complexes ii.
\newblock {\em Ann. Scuola Norm. Sup. Pisa. (2)}, 1:333--354, 1932.

\bibitem[D'A93]{Dangelo}
John~P. D'Angelo.
\newblock {\em Several complex variables and the geometry of real
  hypersurfaces}.
\newblock Studies in Advanced Mathematics. CRC Press, Boca Raton, FL, 1993.

\bibitem[Fer96]{Ferrand2}
Jacqueline Ferrand.
\newblock The action of conformal transformations on a {R}iemannian manifold.
\newblock {\em Math. Ann.}, 304(2):277--291, 1996.

\bibitem[Fra06]{Charlie2}
Charles Frances.
\newblock {A Ferrand-Obata theorem for rank one parabolic geometries}, 2006.
\newblock arXiv:math.DG/0608537.

\bibitem[FS74]{FS}
G.~B. Folland and E.~M. Stein.
\newblock Estimates for the {$\bar \partial \sb{b}$} complex and analysis on
  the {H}eisenberg group.
\newblock {\em Comm. Pure Appl. Math.}, 27:429--522, 1974.

\bibitem[FT02]{Charlie1}
Charles Frances and C{\'e}dric Tarquini.
\newblock Autour du th\'eor\`eme de {F}errand-{O}bata.
\newblock {\em Ann. Global Anal. Geom.}, 21(1):51--62, 2002.

\bibitem[Gam01]{Gam}
Najoua Gamara.
\newblock The {CR} {Y}amabe conjecture---the case {$n=1$}.
\newblock {\em J. Eur. Math. Soc. (JEMS)}, 3(2):105--137, 2001.

\bibitem[GY01]{GY}
Najoua Gamara and Ridha Yacoub.
\newblock C{R} {Y}amabe conjecture---the conformally flat case.
\newblock {\em Pacific J. Math.}, 201(1):121--175, 2001.

\bibitem[Heb97]{Heb}
Emmanuel Hebey.
\newblock {\em Introduction \`{a} l'analyse non lin\'{e}aire sur les
  vari\'{e}t\'{e}s}.
\newblock Fondations. Diderot Sciences, 1997.

\bibitem[Jac90]{Jacobowitz}
Howard Jacobowitz.
\newblock {\em An introduction to {CR} structures}, volume~32 of {\em
  Mathematical Surveys and Monographs}.
\newblock American Mathematical Society, Providence, RI, 1990.

\bibitem[JL87]{JL1}
David Jerison and John~M. Lee.
\newblock The {Y}amabe problem on {CR} manifolds.
\newblock {\em J. Differential Geom.}, 25(2):167--197, 1987.

\bibitem[JL89]{JL2}
David Jerison and John~M. Lee.
\newblock Intrinsic {CR} normal coordinates and the {CR} {Y}amabe problem.
\newblock {\em J. Differential Geom.}, 29(2):303--343, 1989.

\bibitem[Kam93]{Kamishima1}
Yoshinobu Kamishima.
\newblock A rigidity theorem for {CR} manifolds and a refinement of {O}bata and
  {L}elong-{F}errand's result.
\newblock In {\em Geometry and its applications (Yokohama, 1991)}, pages
  73--83. World Sci. Publ., River Edge, NJ, 1993.

\bibitem[Kam96]{Kamishima2}
Yoshinobu Kamishima.
\newblock Geometric flows on compact manifolds and global rigidity.
\newblock {\em Topology}, 35(2):439--450, 1996.

\bibitem[Kob95]{Kobayashi}
Shoshichi Kobayashi.
\newblock {\em Transformation groups in differential geometry}.
\newblock Classics in Mathematics. Springer-Verlag, Berlin, 1995.
\newblock Reprint of the 1972 edition.

\bibitem[Kup94]{Kuperberg}
Krystyna Kuperberg.
\newblock A smooth counterexample to the {S}eifert conjecture.
\newblock {\em Ann. of Math. (2)}, 140(3):723--732, 1994.

\bibitem[Laf88]{Lafontaine}
Jacques Lafontaine.
\newblock The theorem of {L}elong-{F}errand and {O}bata.
\newblock In {\em Conformal geometry (Bonn, 1985/1986)}, Aspects Math., E12,
  pages 93--103. Vieweg, Braunschweig, 1988.

\bibitem[Lee96]{Lee}
John~M. Lee.
\newblock C{R} manifolds with noncompact connected automorphism groups.
\newblock {\em J. Geom. Anal.}, 6(1):79--90, 1996.

\bibitem[LF71]{Ferrand1}
Jacqueline Lelong-Ferrand.
\newblock Transformations conformes et quasi-conformes des vari\'et\'es
  riemanniennes compactes (d\'emonstration de la conjecture de {A}.
  {L}ichnerowicz).
\newblock {\em Acad. Roy. Belg. Cl. Sci. M\'em. Coll. in--8$\deg\ $(2)},
  39(5):44, 1971.

\bibitem[LP87]{LP}
John~M. Lee and Thomas~H. Parker.
\newblock The {Y}amabe problem.
\newblock {\em Bull. Amer. Math. Soc. (N.S.)}, 17(1):37--91, 1987.

\bibitem[MZ51]{Montgomery-Zippin}
Deane Montgomery and Leo Zippin.
\newblock Existence of subgroups isomorphic to the real numbers.
\newblock {\em Ann. of Math. (2)}, 53:298--326, 1951.

\bibitem[Oba72]{Obata}
Morio Obata.
\newblock The conjectures on conformal transformations of {R}iemannian
  manifolds.
\newblock {\em J. Differential Geometry}, 6:247--258, 1971/72.

\bibitem[Pan90]{Pansu}
Pierre Pansu.
\newblock Distances conformes et cohomologie $l^n$.
\newblock {\em Publ. Univ. Pierre et Marie Curie}, 92, 1990.

\bibitem[Sch74]{Schweitzer}
Paul~A. Schweitzer.
\newblock Counterexamples to the {S}eifert conjecture and opening closed leaves
  of foliations.
\newblock {\em Ann. of Math. (2)}, 100:386--400, 1974.

\bibitem[Sch95]{Schoen}
R.~Schoen.
\newblock On the conformal and {CR} automorphism groups.
\newblock {\em Geom. Funct. Anal.}, 5(2):464--481, 1995.

\bibitem[Sha97]{Sharpe}
R.~W. Sharpe.
\newblock {\em Differential geometry}, volume 166 of {\em Graduate Texts in
  Mathematics}.
\newblock Springer-Verlag, New York, 1997.
\newblock Cartan's generalization of Klein's Erlangen program, With a foreword
  by S. S. Chern.

\bibitem[Spi97]{Spiro}
Andrea~F. Spiro.
\newblock Smooth real hypersurfaces in {$\mathbf C\sp n$} with a non-compact
  isotropy group of {CR} transformations.
\newblock {\em Geom. Dedicata}, 67(2):199--221, 1997.

\bibitem[Thu97]{Thurston}
William~P. Thurston.
\newblock {\em Three-dimensional geometry and topology. {V}ol. 1}, volume~35 of
  {\em Princeton Mathematical Series}.
\newblock Princeton University Press, Princeton, NJ, 1997.
\newblock Edited by Silvio Levy.

\bibitem[Web77]{Webster1}
S.~M. Webster.
\newblock On the transformation group of a real hypersurface.
\newblock {\em Trans. Amer. Math. Soc.}, 231(1):179--190, 1977.

\bibitem[Web78]{Webster2}
S.~M. Webster.
\newblock Pseudo-{H}ermitian structures on a real hypersurface.
\newblock {\em J. Differential Geom.}, 13(1):25--41, 1978.

\bibitem[Won77]{Wong}
Bun Wong.
\newblock Characterization of the unit ball in {${\bf C}\sp{n}$} by its
  automorphism group.
\newblock {\em Invent. Math.}, 41(3):253--257, 1977.

\bibitem[Won03]{WongSurvey}
B.~Wong.
\newblock On complex manifolds with noncompact automorphism groups.
\newblock In {\em Explorations in complex and Riemannian geometry}, volume 332
  of {\em Contemp. Math.}, pages 287--304. Amer. Math. Soc., Providence, RI,
  2003.

\end{thebibliography}

\end{document}